\documentclass{amsart}

\usepackage{graphicx}
\usepackage{amsmath}
\usepackage{amsfonts}
\usepackage{amssymb}

\newtheorem{theorem}{Theorem}[section]

\newtheorem{proposition}[theorem]{Proposition}

\theoremstyle{definition}

\theoremstyle{remark}

\newcommand{\be}{\begin{equation}}
\newcommand{\ee}{\end{equation}}
\newcommand{\ben}{\begin{equation*}}
\newcommand{\een}{\end{equation*}}
\newcommand{\beq}{\begin{eqnarray}}
\newcommand{\eeq}{\end{eqnarray}}

\newcommand{\bma}{\begin{pmatrix}}
\newcommand{\ema}{\end{pmatrix}}
\newcommand{\bal}{\begin{align}}
\newcommand{\eal}{\end{align}}
\newcommand{\baln}{\begin{align*}}
\newcommand{\ealn}{\end{align*}}
\newcommand{\bald}{\begin{aligned}}
\newcommand{\eald}{\end{aligned}}

\newcommand{\wt}[1]{\widetilde{#1}}
\newcommand{\goth}[1]{\mathfrak{#1}}

\newcommand{\der}{{\rm d}}
\newcommand{\w}{{\scriptstyle\wedge}}

\newcommand{\real}{\mathbb{R}}

\newcommand{\K}{{\mathcal K}}
\newcommand{\M}{{\mathcal M}}

\renewcommand{\P}{\mathcal{P}}

\newcommand{\Ric}{{\rm Ric}}
\newcommand{\id}{{\rm id}}

\renewcommand{\o}{\goth{o}}
\newcommand{\co}{\goth{co}}

\newcommand{\h}{\goth{h}}
\renewcommand{\a}{\goth{a}}

\newcommand{\wcov}{\widehat{\nabla}}

\begin{document}
\title{On three-dimensional Weyl structures with reduced holonomy} 
 
\author{Micha\l~ Godli\'nski} \address{Instytut Fizyki Teoretycznej,
Uniwersytet Warszawski, ul. Ho\.za 69, Warszawa, Poland}
\email{godlinsk@fuw.edu.pl} 

\author{Pawe\l~ Nurowski} \address{Instytut Fizyki Teoretycznej,
Uniwersytet Warszawski, ul. Ho\.za 69, Warszawa, Poland}
\email{nurowski@fuw.edu.pl} 

\date{\today}

\begin{abstract}
Cartan's list of 3-dimensional Weyl structures with reduced 
holonomy is revisited. We show that the only Einstein-Weyl structures on
this list correspond to the structures generated by the solutions of the 
dKP equation.   
\end{abstract}

\maketitle

\allowdisplaybreaks

\section{Introduction}

\noindent 
In Ref. \cite{Car1} Elie Cartan gives a complete list of 3-dimensional
Weyl geometries with reduced holonomy. Cartan does not study the
Einstein-Weyl equations for the geometries from his list. On the other
hand, in recent years, various authors \cite{DMT,H,JT,W} have been
studying the Einstein-Weyl equations in 3-dimensions, mainly due to their
realtions with the twistor theory and the integrable systems theory. In
particular, Dunajski et al \cite{DMT} characterized all 3-dimensional
Einstein-Weyl spaces which admit covariantly constant weighted vector
field, as being generated by solutions to the dispersionless
Kadomtsev-Petviashvili (dKP) equation. 
Their analysis is very much in the spirit of the reduced holonomy
ideas, since the existence of such a vector field reduces the holonomy
of the considered Weyl geometry. However, it is 
not clear from their analysis, if all the 3-dimensional Weyl geometries
with reduced holonomy may be obtained by means of the assumption
about an existence of a covariantly constant weighted vector
field. Quick inspection of the Cartan list of Ref. \cite{Car1} (look
also at Table \ref{t1} of the present paper) shows, that such an 
assumption is very strong and that it excludes a large class of Weyl
geometries with reduced holonomy. A natural question, if the
geometries from this class may be Einstein, is addressed in the present
paper. Here, we first simplify and rephrase in
the modern language Cartan's classification of 3-diemsnional Weyl
geometries with reduced holonomy. This is done by inspecting all 
possible subalgebras of $\co(2,1)$ and $\co(3)$. Then, by
means of the integration of the first structure equations, we
determine which of
them may appear as the Weyl holonomy algebras. The integration
procedure enables us to give canonical representatives of 
the metric and the Weyl potential for each algebra representing the
reduced holonomy. We also specify the geometric object that reduces
the holonomy. It is either a covariantly constant vector field or a
covariantly constant null direction. This second possibility
corresponds to the class of Weyl geometries present at the Cartan
list, but not considered by Dunajski et al. The last part of our paper
imposes the Einstein condition on all the geometries from the Cartan
list. The result is included in Propositions \ref{pr1} and \ref{pr2}
which strengthen the results of Ref. \cite{DMT} to the following
statement:
\newline
\emph{All 3-dimensional Einstein-Weyl geometries with reduced holonomy
are either flat or are generated by the solutions of the dispersionless
Kadomtsev-Petviashvili equation.}

\section{Weyl structures}
\noindent Weyl structure on a real $n$-dimensional manifold  $\M$ consists of a conformal class of metrics 
$[g]$ of signature $(p,q)$ and a torsion-free covariant derivative $\nabla$, such that for each representative $g$ of $[g]$ 
there exists a 1-form $\nu$ satisfying
\be \label{e10}
    \nabla g=-2\nu\otimes g.
\ee
When $g$ changes by $g\to e^{2\phi}g$, then $\nu$ changes by $\nu\to\nu-\der\phi$ so as to leave \eqref{e10} invariant. 
The class of pairs $[g, \nu]$ considered modulo this gauge, uniquely defines the Weyl structure. If $\nu$ is closed, 
then the Weyl structure can be locally reduced to a metric structure by an appropriate gauge, thus we assume 
$\der\nu\neq 0$. 

Let 
$$CO(p,q)=\real_+\times O(p,q)=\{M\in {\rm Matr}_{n\times n}(\real)\,|\,M^TgM=\lambda g, \quad \lambda\in\real_+\} $$ 
denote the Lie group preserving the conformal class $[g]$. Then $[g]$ defines the bundle 
$$CO(p,q)\to\P\to\M,$$ 
a reduction of the bundle of linear frames on $\M$. Any Weyl structure on $\M$ is alternatively defined
by a linear $\co(p,q)$-valued torsion-free connection on $\P$. 
This enables us to apply all the results of the theory of connections \cite{KN} to this case. In paricular, 
the notion of holonomy is well defined and one can study Weyl structures with reduced holonomy. By means of 
the Reduction Theorem (\cite{KN}, p.83) a Weyl structure has its holonomy reduced to some subgroup $H\subset CO(p,q)$
if and only if the Weyl connection is reducible to an $\h$-valued connection on the holonomy bundle $H\to\P'\to\M$ of $\P$. 

Let $(e_i)$ be a frame on $\M$, such that the dual coframe $(\theta^i)$ is orthonormal for some representative $g$ of $[g]$
i.e. $g=g_{ij}\theta^i\theta^j$ with all the coefficients $g_{ij}$
being constant. By a Weyl connection 1-form $\Gamma$, 
we understand the pullback of the Weyl connection from  $CO(p,q)\to\P\to\M$ to $\M$, through the frame $(e_i)$ considered 
as a section of $\P$. The Weyl connection 
1-forms $\Gamma^i_{~j}$ are uniquely defined by the relations
\begin{eqnarray}
 &\der\theta^i+\Gamma^i_{~j}\w\theta^j=0, \label{e20} \\ 
 &\Gamma_{(ij)}=g_{ij}\,\nu, \quad \text{where} \quad \Gamma_{ij}=g_{jk}\Gamma^k_{~j}.\label{e25}
\end{eqnarray}
A Weyl structure has its holonomy reduced to $H$ if the matrix
$(\Gamma^i_j)$ takes values in the Lie algebra $\h$ of $H$.
Due to
\ben
 \co(p,q)=\real\oplus\o(p,q)
\een
$\Gamma$ decomposes into the $\real$-valued part $\nu$ and the 
$\o(p,q)$-valued part $\wt{\Gamma}$ so that
\be\label{e27}
  \Gamma=\nu\cdot\id+\wt{\Gamma}.
\ee 
Hence the subgroups $H\subset O(p,q)\subset CO(p,q)$ do not appear as holonomy groups for Weyl structures. 

A tensor field of weight $m$ on $\M$ is a tensor object $T$ transforming $T\to e^{m\phi}T$ when 
$g\to e^{2\phi}g$. The weighted covariant derivative of a 
$(k,l)$-tensor field of weight $m$  
\ben
   \wcov T= \nabla  T+m\nu\otimes T
\een
is a $(k+1,l)$-tensor field of weight $m$. If $\wcov T=0$, then $T$ is
said to be covariantly constant or constant, for short. 
A direction $\K$ spanned by a vector field $K$ is said to be constant
when $\nabla_X K\in\K$ for an arbitrary vector field $X$. 
A non-null direction $\K$ is constant iff there is a constant vector
field $K\in\K$ of weight $-1$. The existence of 
a constant null direction is a weaker property then the existence of a constant weighted vector in this direction.

The curvature 2-form $\Omega$, the Ricci tensor $\Ric$ and the Ricci
scalar $R$ of a Weyl structure are defined by
\beq
 &\Omega^i_{~j}=\der\Gamma^i_{~j}+\Gamma^i_{~k}\w\Gamma^{k}_{~j}, \nonumber  \\
 & \Omega^i_{~j}=\tfrac{1}{2}\Omega^i_{~jkl}\theta^k\w\theta^l, \nonumber \\
 & \Ric_{ij}=\Omega^k_{~ikj}, \nonumber \\
 & R=\Ric_{ij}g^{ij}. \nonumber
\eeq
$\Omega$ and $\Ric$ have weights 0 whereas $R$ has weight $-2$.
Einstein -- Weyl (E-W) structures are, by definition, those Weyl
structures for which the symmetric trace-free part of the  
Ricci tensor vanishes 
\ben
 \Ric_{(ij)}-\tfrac{1}{n}R\cdot g_{ij} =0.
\een
Weyl structure is flat i.e. $\Omega=0$ iff it has a (local) representative $(\eta, 0)$, where $\eta$ is the flat metric.

\section{Three-dimensional Weyl structures with reduced holonomy}

\noindent 
In order to find all possible 3-d Weyl structures with reduced
holonomy, we 
integrate equations \eqref{e20} for each subalgebra of 
$\co(2,1)$  or $\co(3)$. These subalgebras are classified in \cite{PW}
up to adjoint transformations. We use this classification in the following. 

We begin with the more complicated Lorentzian case.
Let us choose a coframe $(\theta^1,\theta^2,\theta^3)$  such that
$$g=(\theta^2)^2-2\theta^1\theta^3. $$ 
The algebra $\co(2,1)$ now reads
\ben
 \bma  p+a & b & 0 \\
       c & p & b \\
      0 & c & p-a
 \ema. 
\een
The subalgebras with $p\neq 0$ are the following  
\begin{align*}
 && && &A: && c=0, &&&& \\
 && && &B_q: &&  c=0,\, a=-(q+1)p,\, q\in\real, &&&&\\
 && && &C: &&   c=0,\, b=0, &&&&\\
 && && &D_q: &&   c=0,\, b=0,\, a=(q+1)p,\, q>-1, &&&&\\
 && && &E: && a=0,\, c=-b, &&&&\\
 && && &F: &&  c=0,\, a=0,\, b=\pm\,p, &&&&\\
 && && &G_q: &&  a=0,\, b=qp,\, c=-qp,\, q\in\real. &&&&
\end{align*}
Obviously, $A$ contains $B,C,D,F$, and $E$ contains $G$. 

Let us integrate the system \eqref{e20} for subalgebra $A$. In this
case the system reads
\be\label{e70}
 \bald
 &\der\theta^1+(\nu+\alpha)\w\theta^1+\beta\w\theta^2=0, \\
 &\der\theta^2+\nu\w\theta^2+\beta\w\theta^3=0, \\
 &\der\theta^3+(\nu-\alpha)\w\theta^3=0,
\eald
\ee
where $\alpha=\tfrac{1}{2}(\Gamma^1_{~1}-\Gamma^2_{~2})$, $\beta=\Gamma^1_{~2}=\Gamma^2_{~3}$, $\gamma=\Gamma^2_{~1}=\Gamma^2_{~2}$,
$\nu=\Gamma^2_{~2}=\tfrac{1}{2}(\Gamma^1_{~1}+\Gamma^2_{~2})$. We have a three-parameter family of transformations
preserving this system; this is the Lie group $G_A$ of
algebra $A$. The coframe transformation $\theta^i \to
M^i_{~j}\theta^j$ with $M\in G_A$ of the form 
$M=\exp\left( \begin{smallmatrix} t&0&0\\0&0&0\\0&0&-t \end{smallmatrix} \right )$ sends $\alpha\to \alpha-\der t$, $\beta\to e^t\beta$ 
and $\nu\to\nu$. Similarly $M=\exp\left( \begin{smallmatrix} t&0&0\\0&t&0\\0&0&t \end{smallmatrix} \right )$
sends $\alpha\to \alpha$, $\beta\to\beta$, $\nu\to\nu-\der t$ and 
$M=\exp\left( \begin{smallmatrix} 0&t&0\\0&0&t\\0&0&0 \end{smallmatrix} \right )$ transforms $\beta\to\beta-t\alpha-\der t$, leaving 
$\alpha$ and $\nu$ invariant. Exploiting this gauge freedom we
easily achieve 
\be\label{e80}
 \beta\w\theta^2\w\theta^3=0,\qquad (\nu+\alpha)\w\theta^1\theta^3=0, \qquad (\nu-\alpha)\w\theta^3=0.
\ee
Now $\der\theta^3=0$, which enables us to make $\theta^3=\der
x$. Moreover, $\der\theta^2\w\theta^2=0$ and
$\der\theta^1\w\theta^3=0$, so that $\theta^2=a\der z$
and $\theta^1=\der y +b\der x$. In addition, 
$0\neq \theta^1\w\theta^2\w\theta^3=ab\der x\w\der y\w\der z$, hence
$(x,y,z)$ is a coordinate system on $\M$. In this coordinate system, 
using \eqref{e70}, \eqref{e80}, it is easy to get 
$\nu=-(\log a)_y\der y+c \der x$ and $g=a^2\der z^2-2\der x\der y-2b\der x^2$. This, when rescaled
via $g\to a^{-2}g$, $\nu\to\nu-\frac{1}{2a}\der a$, after an
apropriate redefinition of  $a,b,c$ read
\be \label{e90}
 g=\der z^2+2H(x,y,z)\der x\der y+K(x,y,z)\der x^2, \quad \nu=L(x,y,z)\der x-\tfrac{1}{2}\frac{H_z}{H}\der z,
\ee
where $H,K$ are sufficiently smooth arbitrary functions of the
coordinates $(x,y,z)$. In the above gauge the remaining connection
1-forms $\alpha$, $\beta$ are:
\be\label{e100}
 \bald
 \alpha&= (\tfrac{1}{2}\frac{K_y}{H}-\frac{H_x}{H}-L)\der x-\tfrac{1}{2}\frac{H_z}{H}\der z,\\
 \beta&= \tfrac{1}{2}(\frac{H_z}{H}K-K_z)\der x+L\der z.
 \eald
\ee
\noindent
Since the subgroup $G_A$ of $CO(2,1)$
preserves a null direction, these Weyl structures have a constant null
direction. It is generated by the vector field $\partial_y$.

\begin{table}[t]
\centering
\begin{tabular}{| c | c | c | c |} \hline
  &&&holonomy   \\
  type & structure & constant object & algebra \\ \hline  &&& \\

  & $g=\der z^2+2H(x,y,z)\der x\der y+K(x,y,z)\der x^2$ && \\  $A$ && null direction of $\partial_y$ & $\a_1\oplus\real$ \\
  & $\nu=L(x,y,z)\der x-\frac{1}{2H}H_z\der z$ && \\  &&& \\ \hline &&& \\

  & $g=\der z^2+2\der x\der y+K(x,z)\der x^2$ && \\  $B_0$ && null vector $\partial_y$ & $\a_1$ \\
  & $\nu=L(x,y,z)\der x$ && \\  &&& \\ \hline &&& \\

  & $g=\der z^2+2H(x,y,z)\der x\der y+K(x,y,z)\der x^2$ &null 1-form $\der x$, & \\ $B_{-2}$ &&  null vector $\partial_y$& $\a_1$ \\
  & $\nu=\frac{1}{4H}(K_y-2H_x)\der x-\frac{1}{2H}H_z\der z$ & of weight $-2$& \\  &&& \\ \hline &&& \\

  $B_{q}$ & $g=\der z^2+2\der x\der y+K(x,y,z)\der x^2$ & null vector $\partial_y$& $\a_1$ {\scriptsize for }${\scriptstyle q\neq -1}$ \\ 
   && of weight $q$ & \\ 
  ${\scriptstyle q\neq 0,-2}$ & $\nu=-\tfrac{1}{2q}K_y \der x$ && $\real^2$ {\scriptsize for }${\scriptstyle q= -1}$\\  
   &&& \\ \hline &&& \\

  & $g=\der z^2+2H(x,y,z)\der x\der y$ & & \\ $C$ && spatial direction of $\partial_z$ & $\real^2$ \\  
  & $\nu=-\tfrac{1}{2H}H_z\der z$ && \\  &&& \\ \hline &&&\\

  &$g=\der z^2+2H(y,z)\der x\der y$ & spatial direction of $\partial_z$, & \\ $D$ &&& $\real$ \\  
  & $\nu=-\tfrac{1}{2H}H_z\der z$ & null vector $\partial_x$ & \\ &&& \\ \hline &&&\\

  &$g=K(x,y,z)(\der x^2+\der y^2)\pm\der z^2$ & & \\ $E$ && (timelike) direction & $\real^2$ \\ 
  &$\nu=-\frac{1}{2K}K_z\der z$& of $\partial_z$ & \\  &&& \\ \hline 
\end{tabular}
\newline \caption{Three-dimensional Weyl structures with reduced holonomy} \label{t1}
\end{table}

Let us pass to the structures with holonomy $B_q$. Since $B_q$ is
contained in $A$ we can use (\ref{e90})-(\ref{e100}) with the condition of further reduction of the holonomy.  
This is reduced from $A$ to $ B_q$ iff $\alpha=-(q+1)\nu$, which
restricts the possible $H$, $K$ and $L$ by 
\ben
  (q+2)H_z=0,\qquad \text{and}\qquad qL=\frac{H_x}{H}-\tfrac{1}{2}\frac{K_y}{H}.  
\een
For $q=-2$ we have
\be\label{e110}
 g=\der z^2+2H(x,y,z)\der x\der y+K(x,y,z)\der x^2, \quad \nu=(\frac{K_y}{4H}-\frac{H_x}{2H})\der x-\frac{H_z}{2H}\der z.
\ee
When $q\neq -2$, $H=H(x,y)$ and it may be gauged to $H=1$ by means of the 
transformation $H\to Y_y(x,y)$, $K\to
2Y_x+2K(x,y,z)$ followed by the change of coordinates $y\to Y$. In
this gauge $2qL=-K_y$. Thus, for $q\neq 0,-2$ we have
\ben
 g=\der z^2+2\der x\der y+K(x,y,z)\der x^2, \quad \nu=-\tfrac{1}{2q}K_y \der x,
\een
and for $q=0$ 
\be \label{e120}
 g=\der z^2+2\der x\der y+K(x,z)\der x^2, \quad \nu=L(x,y,z)\der x.
\ee
All structures with holonomy $B_q$ have constant null vector field of
weight $q$. In the above coordinates it is given by $\partial_y$. 
In particular, in \eqref{e120} $q=0$, thus we have a constant null
vector field $\partial_y$ there; in \eqref{e110} $q=-2$ and we have a 
constant null 1-form $\der x$ in this case.   

We find structures with holonomy $C$ and $D$ for $q=0$ in an analogous
way. We show that if the holonomy is reduced to types $D$ for $q\neq
0$ and  $F$, then the corresponding Weyl structures are neccessarily
flat. In the nontrivial cases of Weyl structures with holonomies of
type $C$ and $D$ with $q=0$ we have a constant spatial direction. The
case $D$ with $q=0$ admits also a constant null vector. In a similar way, 
we get a family of Weyl structures with holonomy of type $E$
\be\label{e150}
 g=K(x,y,z)(\der x^2+\der y^2)-\der z^2, \quad  \nu=-\tfrac{K_z}{2K}\der z.
\ee
They admit a constant timelike direction generated by $\partial_z$. We
close the discussion of the Lortentzian case by mentioning that the 
structures with holonomy of type $G$ do not exist.

The Euclidean case is much simpler due to the structure of $\co(3)$. It has only two proper subalgebras
up to adjoint automorphisms. They constitute the counterparts of types
$E$ and $G$ from the Lorentzian case. Structures of type $G$ do not
exist, and structures of type $E$ have the form similar to
\eqref{e150}, differing from it merely by the sign standing by 
the $\der z^2$ term.   

All the structures with reduced holonomy, together with their
geometric characterization
are given in Table \ref{t1}. Types $A$ -- $D$ have
Lorentzian signature, type $E$ may have both the Lorentzian and the
Euclidean signatures. In this table $\a_1$
denotes the unique 2-dimensional non-comutative Lie algebra.

\subsection{Three-dimensional E-W structures with reduced holonomy}
We calculate E-W equations in three dimensions 
\ben
  \Ric_{(ij)}-\tfrac{1}{3}R\cdot  g_{ij}=0
\een
for the structures of Table \ref{t1}. It appears, as it was observed
in \cite{DMT}, that E-W structures of types $B_q$ for $q\neq-\tfrac{1}{2}$, 
$C$, $D$ and $E$ are flat ($\Omega=0$). Type $B_{-1/2}$ case is more
interesting. Here, the E-W equations reduce to the dispersionless 
Kadomtsev-Petviashvili (dKP) equation
\ben
 (KK_y-2K_x)_y=K_{zz}.
\een

The structures of type $A$ were not considered by \cite{DMT}. The E-W
system for them  consists of four PDEs for the functions $H,K,L$. One
of these equations is
$H_{yz}H-H_zH_y=0$ with the general solution $H=H_1(x,z)H_2(x,y)$. We
absorb $H_2(x,y)$ by a redefinition $y=y(x,Y)$, $H_2y_Y=1$, 
of $y$-coordinate. Hence, without loss of generality, we take
$H=H(x,z)$. After a substitution $H=\exp(-F(x,z))$, 
$K=G(x,y,z)\exp(-F(x,z))$, two of the remaining three E-W equations read
$$
 L_y=G_{yy}, \qquad L_z=G_{yz}+\tfrac{1}{2}F_{xz}.
$$
They can be easily solved. Now, the Weyl structure reads
$$ 
 g=e^F\der z^2-2\der x\der y+G\der x^2, \qquad \nu=(G_y+f'(x))\der x.
$$
It turns out, that this structure admits null constant vector field
$X=\exp(\tfrac{3}{4}F-\tfrac{1}{2}f(x))\partial_y$ of weight $-1/2$, so the holonomy is of type 
$B_{-1/2}$. Hence, if we impose the last of the E-W equations, the structure will reduce to the one 
generated by the solutions of the dKP equation. 
Thus, type $A$, although more general than $B_{-1/2}$, provides no
essential generalization of the dKP equation. 
We may summarize this section with the following two propositions.
\begin{proposition}\label{pr1}
 Every three-dimensional Euclidean Einstein-Weyl geometry with reduced holonomy is flat.
\end{proposition} 
\begin{proposition}\label{pr2}
 Every three-dimensional Lorentzian Einstein-Weyl geometry with reduced holonomy is flat or has a constant null vector field of weight 
 $-\tfrac{1}{2}$. In the latter case E-W equations reduce to the dKP equation in some coordinate system.
\end{proposition}

\end{document}